\documentclass[letterpaper, 10 pt, conference]{ieeeconf}  
\IEEEoverridecommandlockouts                  
\usepackage{amsthm}
\usepackage{cite}
\usepackage{amsmath,amssymb,amsfonts}
\usepackage{algorithmic}
\usepackage{graphicx}
\usepackage{textcomp}
\usepackage{xcolor}
\usepackage{defs_new}
\usepackage{caption}
\usepackage{subcaption}
\usepackage{todonotes}
\usepackage{multicol}
\usepackage{mathtools}
\usepackage{todonotes}
\usepackage{float}
\usepackage{booktabs}
\usepackage{array, multirow}
\usepackage[numbered,framed]{matlab-prettifier}
\usepackage{algorithm} 
\usepackage{algorithmic}  
\usepackage[linesnumbered,ruled,lined,algo2e]{algorithm2e}
\usepackage{tabularray}

\usepackage{tikz}
\usetikzlibrary{arrows.meta}
\usetikzlibrary{positioning, calc}
\usetikzlibrary{shadows, fit, shapes.multipart, shapes.callouts}

\usetikzlibrary{decorations.pathmorphing}

\usepackage[framemethod=tikz,xcolor=true]{mdframed}

\definecolor{metabluee}{HTML}{B5CDF5}
\definecolor{metablue}{HTML}{0064E0}
\definecolor{metafg}{HTML}{1C2B33}
\definecolor{metabg}{HTML}{F1F4F7}
\definecolor{metabgdeep}{HTML}{D9EFFF}
\definecolor{metagreen}{HTML}{EAFFE8}
\definecolor{metagreen}{HTML}{FCFFEE}
\definecolor{metared}{HTML}{FFEAE8}

\definecolor{Orchid}{RGB}{218,112,214}
\definecolor{Bittersweet}{RGB}{254,111,94}
\definecolor{CornflowerBlue}{RGB}{100,149,237}
\definecolor{Goldenrod}{RGB}{218,165,32}
\definecolor{LimeGreen}{RGB}{50,205,50}
\definecolor{Teal}{RGB}{0,128,128}

\newmdenv[backgroundcolor=metabgdeep, roundcorner=10pt, skipabove=4pt, linewidth=0pt, innertopmargin=4pt]{myframe}

\newmdenv[backgroundcolor=metabg, roundcorner=6pt, skipabove=2pt, linewidth=0pt, innertopmargin=0pt]{myOCP}

\newmdenv[backgroundcolor=metared, roundcorner=10pt, skipabove=7pt, linewidth=0pt, innertopmargin=7pt]{myalgo}

\usepackage{tikz}
\def\BibTeX{{\rm B\kern-.05em{\sc i\kern-.025em b}\kern-.08em
    T\kern-.1667em\lower.7ex\hbox{E}\kern-.125emX}}

\title{\LARGE \bf
Data-Driven Density Steering via the Gromov-Wasserstein Optimal Transport Distance}

\author{Haruto Nakashima, Siddhartha Ganguly, Kenji Kashima 
\thanks{The authors are with the Applied Mathematics and Physics Course, Graduate School of
Informatics, Kyoto University, Kyoto 606-8501, Japan. Email: nakashima.haruto.72v@st.kyoto-u.ac.jp, ganguly.siddhartha.7p@kyoto-u.ac.jp, kk@i.kyoto-u.ac.jp. The authors acknowledge the JSPS KAKENHI under Grant Number JP21H04875.%
}
}

\begin{document}

\maketitle
\thispagestyle{empty}
\pagestyle{empty}

\begin{abstract}
We tackle the data-driven chance-constrained density steering problem using the Gromov-Wasserstein metric. The underlying dynamical system is an unknown linear controlled recursion, with the assumption that sufficiently rich input-output data from pre-operational experiments are available. The initial state is modeled as a Gaussian mixture, while the terminal state is required to match a specified Gaussian distribution. We reformulate the resulting optimal control problem as a difference-of-convex program and show that it can be efficiently and tractably solved using the DC algorithm. Numerical results validate our approach through various data-driven schemes. 
\end{abstract}


\section{Introduction}\label{sec:intro}

The term \emph{data-driven} has become increasingly prevalent in the modern control literature \cite{ref:DDCont:main:CDP:PT}. In numerous control applications, the underlying dynamical system and its governing physics are often unknown, with the designer having access solely to input-output data. In such contexts, control designs can be realized through two primary methodologies: an \emph{indirect} approach, where a system model is first identified prior to the synthesis of a controller, and a \emph{direct} approach, which circumvents model identification altogether, generating control commands directly from the available data. 

Within the framework of the latter data-driven paradigm, we focus on a class of optimal control problems known as density/covariance steering problems. Density steering aims to simultaneously control the mean and covariance of a dynamical system's state distribution so that they align with a specified target distribution. The steering problem, a decade-old design framework, has witnessed a surge in research interest in recent years, driven by its critical role in applications demanding precise uncertainty quantification and robust performance such as in \cite{ref:DZ:JR:PT:AA:BSP:21}, autonomous driving \cite{ref:JK:PT:CovUnknParam:2023}, and UAV control \cite{ref:GR:PT:StocCont:UAV}. Although covariance steering enables the controller to be designed offline and typically does not require heavy online computation, traditional methods depend on accurate system models --- a requirement that is often difficult to meet in practice. Data-driven techniques offer an attractive alternative by leveraging experimental data to infer system dynamics and design appropriate data-based controllers. We note that, integration of covariance steering with data-driven control has been attempted in the literature \cite{ref:DDCovSteer:JP:PT,ref:L4DC:DDRobCovSteer:JP:PT}. In the present work, we focus on a data-driven formulation of the covariance steering problem, where the terminal cost is defined via the \emph{Gromov-Wasserstein} (GW) distance. This optimal transport metric provides a flexible means of quantifying discrepancies between probability distributions, making it particularly suitable for comparing the intrinsic geometries of the target and achieved state distributions. Importantly, the GW distance is \emph{invariant} to \emph{isometric Euclidean} transformations, such as rotations and translations, which makes it highly effective in scenarios where the relative shape of the distribution is more significant than its absolute position or orientation. 


To motivate the present work, consider the state distribution as an ensemble of particles or a multi-robotic swarm. In many such applications, achieving a specific formation is essential, although its location and orientation are often immaterial. For instance, a formation may require the agents to align in a single row within a two-dimensional region. While previous formulations can realize a desired configuration with a fixed orientation and position, they do not address optimization with respect to rotation. To overcome this limitation, we introduce a novel density control problem that incorporates the Gromov-Wasserstein distance as the terminal cost. By leveraging the invariance property of the GW distance, our formulation is able to simultaneously optimize the dynamical steering and the rotation of the target shape, standing in clear contrast to existing methods. We note that the GW distance has been utilized in the context of density and formation control for known linear systems without chance constraints; see \cite{ref:KM:KK:MinDenLin:2024, ref:NH:SG:KM:KK:Formation:shape:GW}.

\subsection*{Our contributions}

\begin{itemize}[label=\(\circ\), leftmargin=*]

\item In this article we consider the problem of density steering via the Gromov-Wasserstein terminal cost for an unknown linear system in the data-driven setting. In our setting the initial state is given by a Gaussian Mixture Model (GMM) \cite{ref:IMB:BE:GMM:2024, ref:KumKen:CCGMM:CDC} while the target density follows a Gaussian distribution. 

\item To induce tractability we enforce an affine feedback parametrization of the control policies. Under this parametrization, we show that the state remains distributed as a GMM. Our setting also involves chance constraints on the state variable which we reformulate in a deterministic way employing standard risk allocation techniques \cite{ref:OnoWilliams:Risk:Allocation}. 

\item While computing the GW distance for arbitrary distributions is a well-known NP-complete problem, recent work has demonstrated that the Gaussian Gromov-Wasserstein (GGW) distance --- a relaxation of the GW distance for normal distributions --- admits a closed-form expression \cite{ref:GGW:delon2022gromov}. Leveraging this result along with the data-driven parametrization, we show that the ensuing density steering problem can be translated into a difference-of-convex (DC) programming problem. This problem is solved using the DC algorithm (DCA), which iteratively approximates the problem via convex relaxation. Notably, the convexified formulation can be transformed into a semidefinite programming (SDP) problem, allowing it to be solved efficiently using standard convex programming solvers.

\item To streamline our investigation, we consider the case of noiseless data and demonstrate the applicability of the proposed approach to noise-corrupted data solely through numerical simulations. We note that, while the current paper emphasizes practical implementation and simulation results, a rigorous mathematical treatment of the underlying theoretical aspects with corrupted data is beyond the scope of this work and is left for future research.

\end{itemize}
The article unfolds as follows: Section \S\ref{sec:GW_description} lays the foundation by introducing the concept of GW distance, while \S\ref{sec:problem} formulates the minimum-energy GW density steering problem. In \S\ref{sec:main_results}, we unveil our key findings, establishing distributional properties of state evolution and proving that the density steering problem naturally fits within a DC programming framework. Finally, we validate and illustrate our results via a numerical example in \S\ref{sec:numerics}.

\textbf{Notations:} Our notations are standard. \(\mathbb{N}\) and \(\Z\) denotes the natural numbers and the integers, respectively; \(\N \Let \aset{1, 2, \ldots}\) denote the set of positive integers. Let \(\hat{a},\ol{a} \in \N\) and \(\ol{a}>\hat{a}\), then \([\hat{a};\ol{a}]\Let \aset[]{\hat{a},\hat{a}+1,\ldots,\ol{a}}\). Given \(\stdim \in \N\), the set of \((\stdim \times \stdim)\) positive definite and semi-definite matrices are denoted by \(\pdmat{\stdim}\) and \(\psdmat{\stdim}\); for such matrices \(\norm{\cdot}_F\) denotes the standard Frobenius norm. Expectation operator will be denoted by \(\expect(\cdot)\); the probability density function (PDF) of \(y \in \Rbb^{\stdim}\) at \(\hat{y} \in \Rbb^{\stdim}\) is given by \(\prob_y(\hat{y})\). For a Gaussian-distributed random variable \(x \in \Rbb^{\stdim}\) with mean \(\mean \in \Rbb^{\stdim}\) and covariance \(\covar \in \Rbb^{\stdim \times \stdim}\) we will write \(x\sim \mathcal{N}_d(\mean,\covar)\) and in this case its PDF will be denoted by \(\prob_{\mathcal{N}}(\hat{y};\mean,\covar)\). For two random variables \(\hat{x} \in \Rbb^{\stdim}\) and \(\check{y} \in \Rbb^{\dimcon}\), given \(y = \check{y}\) the conditional random variable \(\hat{y}\) will be denoted by \(\hat{y}|y =\check{y}\). We denote the standard probability simplex by \(\simplex{d} \Let \aset[]{\zeta \in \Rbb^{\stdim} \suchthat \zeta_i \ge 0,\, \sum_{i=1}^{\stdim} \zeta_i = 1}\). For any rectangular matrix \(A \in \Rbb^{d_1 \times d_2}\), the Moore-Penrose pseudoinverse is denoted by \(A^\dagger\). For any \(x \in \Rbb^{d}\) and \(M \in \Rbb^{d \times d}\), we denote the standard quadratic form \(\inprod{x}{Mx}\) by \(\norm{x}^2_{M}\).
\section{Background on the Gromov-Wasserstein distance}\label{sec:GW_description}
Let \(\measX(\cdot)\) and \(\measY(\cdot)\) be two given probability density functions specified on metric spaces \(\metspaceX\) and \(\metspaceY\) with respective metrics \(\costX: \metspaceX \times \metspaceX \lra \Rbb\) and \(\costY:\metspaceY \times \metspaceY \lra \Rbb\). Define the set of probability measures with marginals \(\measX(\cdot)\) and \(\measY(\cdot)\) by
\[
\probmeas(\measX,\measY)\Let \left\{\marginal(x,y)\;\middle\vert\;  
\begin{array}{@{}l@{}}
\int_{\metspaceY} \marginal(x,y)\,\dd y = \measX(x) \text{ and}\\ \int_{\metspaceX}\marginal(x,y)\,\dd x = \measY(y)
\end{array}
\right\}.    
\]
Let \(\mathsf{F}\Let \metspaceX^2 \times \metspaceY^2\), and define \(\widehat{d}(x,x',y,y')\Let \abs{\costX(x,x') - \costY(y,y')}^{p}\); then for \(p\ge 1\), the GW distance, \(\GW_{p} \bigl( \measX, \measY\bigr)^{p}\) is defined by the variational problem
\begin{myOCP}
    \begin{align}\label{eq:cont_GW_dist}
   \inf_{\marginal(\cdot) \in \probmeas(\measX,\measY)} \int_{\mathsf{F}} \widehat{d}(x,x',y,y') \dd\marginal(x,y)\dd\marginal(x',y').
\end{align}
\end{myOCP}
We refer the readers to \cite{ref:FM:GW:2011} for an in-depth analysis of the geometric and metric properties of \(\GW(\cdot,\cdot)\). Another standard OT-distance is the Wasserstein distance \(W_p(\measX, \measY)\) that measures the optimal transport cost between two probability measures \(\measX\) and \(\measY\) defined on a common metric space \((\metspaceX,d)\), which is given by  
\begin{align}\label{eq:cont_W_dist}
\mathsf{W}_p(\measX, \measY)^p \Let  \inf_{\marginal \in \probmeas(\measX,\measY)} \int_{\metspaceX^2} d(x,y)^p \;\dd\marginal(x,y).
\end{align}
It is clear from \eqref{eq:cont_GW_dist} and \eqref{eq:cont_W_dist} that while the Wasserstein distance compares distributions on the same metric space by directly transporting mass between points, the GW distance compares distributions on different spaces by aligning their intrinsic geometric structures. GW is invariant to isometric transformations, unlike the Wasserstein distance and it is suited for comparing structured data such as graphs and shapes. We will use a specialized variant of the GW distance, which is detailed below.

\subsection*{The Gaussian GW distance} The optimization problem \eqref{eq:cont_GW_dist} is a highly non-convex mathematical program and is generally extremely difficult to optimize. Even in the discrete case, when \( \measX \) and \( \measY \) are replaced by discrete probability measures \( \measX \Let \sum_{i=1}^n \sourcevec_i \delta_{x_i} \) and \( \measY \Let \sum_{j=1}^{m} \tarvec_j \delta_{y_j} \) where \( \sourcevec \Let (\sourcevec_1, \ldots, \sourcevec_n)^{\top} \in \simplex{n} \) and \( \tarvec \Let (\tarvec_1, \ldots, \tarvec_m)^{\top} \in \simplex{m} \), the resulting problem translates to a quadratic assignment problem (QAP) (if the costs are quadratic, and the transport plan is given by a permutation matrix), which retains the same NP-complete non-convex structure. However, when the distributions are Gaussian, a closed-form expression can be obtained.

\begin{myOCP}
\begin{theorem}\cite[\S 4, Theorem 4.1]{ref:GGW:delon2022gromov}
Let \(d_1,d_2 \in \N\) and without any loss of generality suppose that \(d_1 \ge d_2\). Let \(\measX \Let   \mathcal{N}_{d_1}(\mean_1,\covar_1)\) and \(\measY \Let  \mathcal{N}_{d_2}(\mean_2,\covar_2)\) be two Gaussian measures on \(\Rbb^{d_1}\) and \(\Rbb^{d_2}\). Let \(P_0,D_0\) and \(P_1,D_1\) be the respective diagonalizations of \(\covar_1 (= P_0 D_0 P_0^{\top})\) and \(\covar_1 (= P_1 D_1 P_1^{\top})\) which sort eigenvalues in decreasing order. We suppose that \(\covar_1\) is non-singular (\(\measX\) is not degenerate). Let \(\mathcal{N}_{d_1+d_2}\) be the set of Gaussian measures on \(\Rbb^{d_1+d_2}\), define \(\mathcal{S}\Let \probmeas(\measX,\measY) \cap \mathcal{N}_{d_1+d_2}\), and let \(\widehat{d}(x,x',y,y')\Let \bigl(\norm{x-x'}^2_{\Rbb^{d_1}} - \norm{y-y'^2}_{\Rbb^{d_2}}\bigr)^2\). Then the optimization problem
\begin{align}\label{eq:cont_GGW_dist}
  \hspace{-2mm}\GGW_2^2(\measX,\measY) \Let  \hspace{-2mm}\inf_{\marginal(\cdot) \in \mathcal{S}} \hspace{-1mm}\int_{\mathsf{F}} &\widehat{d}(x,x',y,y')\dd\marginal(x,y)\dd\marginal(x',y')\nn
\end{align}
admits a solution and 
\begin{align}
 \hspace{-2mm}  &\GGW_2^2(\measX,\measY) =  \hspace{1mm}4 \bigl( \tr(D_0)- \tr (D_1)\bigr)^2 \nn \\& + 8 \norm{D_0^{(n)} - D_1}^2_{F}  - 8 \bigl(\norm{D_0^{(n)}}^2_{F} -  \norm{D_1}^2_{F}\bigr).
\end{align}
\end{theorem}  
\end{myOCP}

\section{Problem setup}\label{sec:problem}
We now formulate our data-driven GGW density steering problem. Consider the time-invariant discrete-time linear dynamical systems
\begin{align}\label{LTI_system}
    x_{t+1} = A x_t + B u_t  \quad  \text{for }t\in\Nz
\end{align}
with the following data:
\begin{enumerate}[label=\textup{(\alph*)}, leftmargin=*, widest=b, align=left]
\item \label{ocpdata:1} \(x_t \in \Rbb^{\stdim}\) and \(u_t \in \Rbb^{\condim}\) denotes the state and control vectors at time \(t \in \Nz\); \(A\in\Rbb^{{\stdim}\times {\stdim}}\) is the system matrix and \(B\in\Rbb^{{\stdim}\times \condim}\) is the actuation matrix, which are assumed to be constant, but \emph{unknown}. The pair \((A,B)\) is controllable. 
\item \label{ocpdata:2} as standard in the data-driven setting, we assume that the \emph{persistently exciting} input-output data \(\dataset \Let \aset[]{x_t^{\dataindex},u_t^{\dataindex},x_{\horizon}^{\dataindex}}_{t=0}^{T-1}\) is available, which is typically collected from offline experiments conducted on the dynamical system. The uncertainty resides in the initial state \(x_0\) and given by the \emph{Gaussian Mixture Model} (GMM). More precisely, let \(\mathsf{K} \in \Nz\); consider the finite sequence of weights \((\alpha_i)_{i=0}^{\mathsf{K}-1}\in\Delta_{\mathsf{K}}\), and covariance matrices \(\bigl(\covar_0^i\bigr)_{i=0}^{\mathsf{K}-1}\) we have 
\begin{align}
    x_0 \sim \GMM \bigl(\{\alpha_i,\mean_0^i,\covar_0^i\}_{i=0}^{\mathsf{K}-1} \bigr).
\end{align}
Moreover, the final state \(x_{\horizon}\) admits a Gaussian distribution with \(\covar_{\horizon}\), where \(\covar_{\horizon} \in \pdmat{d}\).

\item \label{ocpdata:3} Let \(j \in \aset[]{1,\ldots,L}\) and \((a_j,b_j) \in \Rbb^{\stdim} \times \Rbb\) for each \(j \in \aset[]{1,\ldots,L}\). Define the half-plane
\begin{align}\label{eq:halfplane-def}
     H_j \Let \aset[]{x\in\Rbb^{\stdim} \suchthat a_j^{\top} x \le b_j}.
 \end{align}
We enforce half-plane type of chance constraints on the state variable
\begin{align}\label{eq:chance_cons}
\prob\biggl(\bigwedge_{j=1}^{L}\bigwedge_{t=1}^{N}x_t \in H_j\biggr)\ge 1 - V.
\end{align}
Without loss of generality, we assume that \(\norm{a_j} = 1\) in \ref{eq:halfplane-def}, and \(V\) is the allowed probability of violation for each constraint and we have \(0 \le V \le 0.5\). This assumption is justified, as most chance-constrained problems typically permit a risk level below 0.5; see \cite{ref:OnoWilliams:Risk:Allocation, ref:KumKen:CCGMM:CDC}.

\end{enumerate}
\vspace{1mm}
    We wish to steer the terminal density \(\st_{\horizon}\) so that its covariance is isometric to that of the target Gaussian \(x_d \sim \mathcal{N}(0,\covar_d)\), where \(\covar_d \in \psdmat{d}\), while ensuring that all the data \ref{ocpdata:1}--\ref{ocpdata:3} are satisfied. For \(i \in [1;\horizon-1]\) define the sequence of maps \(\pi(\cdot) \Let \pi_{i}(\cdot)\) and let \(U \Let (u_0,\ldots,u_{\horizon-1})^{\top}\). Under the problem data \ref{ocpdata:1}--\ref{ocpdata:3}, we consider the optimal control problem
    \begin{myOCP}
        \begin{equation}
	\label{eq:OCP}
\begin{aligned}
&  \hspace{-2mm}\inf_{\pi(\cdot)}	&&   \hspace{-2mm}J\bigl(U\bigr)\Let \expect \biggl[ \sum_{t=0}^{\horizon-1} \norm{u_t}^2_{R_t} \biggr] + \eps\, \GGW_2^2 (\covar_d,\covar_N) \\
&  \hspace{-2mm} \sbjto		&&  \hspace{-2mm}\begin{cases}
\text{dynamics }\eqref{LTI_system},\, x_0 \sim \GMM\bigl(\{\alpha_i,\mu_0^i,\Sigma_0^i\}_{i=0}^{K-1}\bigr),\\
u_t = \pi_t(x_t), \text{ and }
\text{chance constraints }\eqref{eq:chance_cons},
\end{cases}
\end{aligned}
\end{equation}
    \end{myOCP}
where for all \(t \in [0;\horizon-1]\) the matrix \(R_t \in \pdmat{m}\) and \(\eps>0\) is a parameter.


\section{Main results: data-driven formulation}\label{sec:main_results}
\label{sec:main_result}
The OCP \eqref{eq:OCP} in its current form is numerically intractable as the optimization is infinite-dimensional; indeed, the decision variable \(\pi(\cdot)\) belongs to some class of functions. To induce tractability, following \cite{ref:KumKen:CCGMM:CDC} and \cite{ref:IMB:BE:GMM:2024}, we choose the affine feedback control policy
\begin{align}\label{eq:control:param}
    &\pi_t(x_t) = v_t^i + K_t^i(x_t - \mean_t^i) \quad\text{w.p. }\gamma_i(x_t) \\
    &\text{where} \quad \gamma_i(x_t) \Let \frac{\alpha_i \prob_\mathcal{N}(x_t;\mean_t^i,\covar_t^i)}{\sum_{i=0}^{\mathsf{K}-1} \alpha_i \prob_\mathcal{N}(x_t;\mean_t^i,\covar_t^i)},\nn
\end{align} 
where $v_t^i\in\Rbb^\condim$ and $ K_t^i \in \Rbb^{\condim \times \stdim}$ are probabilistic feedback and feedforward term. 
\subsection{The (state, action)-pair are \(\GMM\)-distributed}
Our control law \eqref{eq:control:param} is inspired by approaches, such as those proposed in \cite{ref:ACC:balci2024density, ref:KumKen:CCGMM:CDC}, but differs slightly in its formulation. First, we verify that under the policy \eqref{eq:control:param}, the state and control variables follow a GMM distribution. The following result is standard; see \cite{ref:KumKen:CCGMM:CDC}.
\begin{lemm}\label{lem:delta-gaussian-integral}
Let \(C \in\Rbb^{p\times m}\) and \(b\in\Rbb^p\). Suppose that \(\delta(\cdot)\) is the standard Dirac distribution. Then for any \(x' \in \Rbb^p\) the following equality holds
\begin{align*}
        \int_{\Rbb^d}\prob_\mathcal{N}(x;\mean,\covar)\delta(x' - (Cx+b))\;\dd x  = \prob_\mathcal{N}(x';\bar{\mean}, \bar{\covar}),
    \end{align*}
where \(\bar{\mean} \Let C\mean+b\) and \(\bar{\covar}\Let C \covar C^{\top}\).
\end{lemm}

\begin{myOCP}
\begin{prop}\label{prop:dist_of_x}
Given \(x_0 \sim \GMM \bigl(\{\alpha_i,\mean_0^i,\covar_0^i\}_{i=0}^{\mathsf{K}-1} \bigr)\) we have \(x_t \sim \GMM(\{\alpha_i,\mean_t^i,\covar_t^i\}_{i=0}^{\mathsf{K}-1})\) for all \(t\in[0;\horizon-1]\), where \(\mean_t^i \in \Rbb^{\stdim}\) and \(\covar_t^i\in \Rbb^{\stdim\times \stdim}\) admit dynamics
\begin{align}
\mean_{t+1}^i &= A \mean_t^i + B v_k^i \label{eq:mean_dynamics}\\
\covar_{t+1}^i &= (A+BK^i_t)\covar_t^i (A+BK^i_t)^\top. \label{eq:covaraicne_dynamics}
\end{align}
\end{prop}  
\end{myOCP}

\begin{proof}
We provide a proof by induction. Assume that \(x_t \sim \GMM\bigl(\{\alpha_i,\mean_t^i,\covar_t^i\}_{i=0}^{\mathsf{K}-1}\bigr)\). Employing the definition of conditional probability densities, the PDF of \(\st_{t+1}\) can be written as follows
\begin{align}
        \prob_{x_{t+1}}(x_{t+1}) = & \int_{\Rbb^{\stdim} \times \Rbb^{\condim}} \prob_{x_{t+1}|x_t = \hat{x}_t,u_t = \hat{u}_t}(x_{t+1}) \nonumber\\
        &\prob_{u_t|x_t = \hat{x}_t}(\hat{u}_t)\prob_{x_t}(\hat{x}_t) \;\dd\hat{u}_t \;\dd\hat{x}_t \label{eq:pdf-of-state}.
\end{align}
In \eqref{eq:pdf-of-state} various terms are given by
\begin{align}
        & \prob_{x_{t+1}|x_t = \hat{x}_t,u_t = \hat{u}_t}(x_{t+1}) = \delta(x_{t+1} - (A\hat{x}_t + B \hat{u}_t)); \nn\\
        & \prob_{u_t|x_t = \hat{x}_t}(\hat{u}_t) = \sum_{i=0}^{\mathsf{K}-1} \gamma_i(\hat{x}_t)\delta(\hat{u}_t - (K_t^i(\hat{x}_t - \mean_t^i)+v_t^i)); \nn \\
        & \prob_{x_t}(\hat{x}_t) = \sum_{i=0}^{\mathsf{K}-1}\alpha_i \, \prob_\mathcal{N}(\hat{x}_t;\mean_t^i,\covar_t^i).\nn
    \end{align}
In \eqref{eq:pdf-of-state}, using the standard evaluation properties of the Dirac delta distribution, we begin by analyzing the integral with respect to \(\hat{u}_t\)
 \begin{align}
        &\int_{\Rbb^\condim} \delta(x_{t+1} - (A\hat{x}_t + B \hat{u}_t))\nn \\  &\times \sum_{i=0}^{\mathsf{K}-1} \gamma_i(\hat{x}_t)\delta(\hat{u}_t - (K_t^i(\hat{x}_t - \mean_t^i)+v_t^i)) \;\dd\hat{u}_t\nn \\ & =\sum_{i=0}^{\mathsf{K}-1}  \gamma_i (\hat{x}_t)\delta(x_{t+1} - ((A + B K_t^i)\hat{x}_t - B K_t^i \mean_t^i + B v_t^i )).\nn
    \end{align}
Plugging this into \ref{eq:pdf-of-state} the assertion follows immediately by invoking Lemma \ref{lem:delta-gaussian-integral}. To wit, 
    \begin{align}
        & \int_{\Rbb^d} \sum_{i=0}^{\mathsf{K}-1} \gamma_i(\hat{x}_t) \delta(x_{t+1} - ((A+BK^i_t)\hat{x}_t - B K_t^i\mean_t^i + Bv_t^i)) \nn\\ & \times \sum_{i=0}^{\mathsf{K}-1} \alpha_i \prob_\mathcal{N}(\hat{x}_t;\mean_t^i,\covar_t^i)d\hat{x}_t = \sum_{i=0}^{\mathsf{K}-1}\alpha_i\int_{\Rbb^\stdim} \prob_\mathcal{N}(\hat{x}_t;\mean_t^i,\covar_t^i) \nn \\   &\delta(x_{t+1} - ((A+BK^i_t)\hat{x}_t - B K_t^i\mean_t^i + Bv_t^i)) \;\dd\hat{x}_t \nn \\
        & = \sum_{i=0}^{\mathsf{K}-1} \alpha_i \prob_\mathcal{N}\bigl(x_{t+1};\check{\mean}, \check{\covar}\big). \nn
    \end{align}
where \(\check{\mean}\Let A \mean_t^i + B v_t^i\) and \(\check{\covar}\Let(A+BK^i_t)\covar_t^i(A+BK^i_t)^\top\). The proof is complete. 
\end{proof}

\begin{myOCP}
\begin{prop}
    Under the policy \eqref{eq:control:param}, the control input is also GMM distributed, i.e.,
    \begin{align}
        u_t \sim \GMM(\{\alpha_i,v_t^i,K_k^i \covar_k^i (K_k^i)^\top\})
    \end{align}
\end{prop}    
\end{myOCP}

\begin{proof}
Utilizing the definition of conditional probability along with the properties of the Dirac distribution, we have
\begin{align*}
    &\prob(u_t)= \int_{\Rbb^\stdim} \prob_{u_t|x_t = \hat{x}_t}(u_t)\prob_{x_t}(\hat{x}_t) \; \dd\hat{x}_t \\
     & = \int_{\Rbb^\stdim} \sum_{i=0}^{\mathsf{K}-1}\gamma_i(\hat{x}_t)\delta(u_t - (K_t^i(\hat{x}_t - \mean_t^i)+v_t^i)) \\
     & \times \sum_{i=0}^{\mathsf{K}-1} \alpha_i \, \prob_\mathcal{N}(\hat{x}_t;\mean_t^i,\covar_t^i) \; \dd\hat{x}_t \\
     & = \sum_{i=0}^{\mathsf{K}-1} \alpha_i \hspace{-1mm}\int_{\Rbb^\stdim} \hspace{-1mm}\delta(u_t - (K_t^i(\hat{x}_t - \mean_t^i)+v_t^i)) \prob_\mathcal{N}(\hat{x}_t;\mean_t^i,\covar_t^i)\; \dd\hat{x}_t \\
     & = \sum_{i=0}^{\mathsf{K}-1}\alpha_i\prob_\mathcal{N}(u_t;v_t^i,K_t^i \covar_t^i (K_t^i)^\top).
\end{align*}
The third equality holds from the definition of $\gamma_{i}(\cdot)$. The last equality follows from Lemma \ref{lem:delta-gaussian-integral}.
\end{proof}
\subsection{A reformulating of the objective function}
Leveraging the parametrization \eqref{eq:control:param}, the function \(J(\cdot)\) in \eqref{eq:OCP} can be decomposed in a \emph{mean} part \(J_{\text{mean}}(\mean_t^i,v_t^i)\) and a \emph{covariance} part \(J_{\text{cov}}(\covar_t^i, K_t^i)\), which are given by 
\begin{align}\label{eq:covariance_objective}
J_{\text{mean}}(\mean_t^i,v_t^i)
& \Let  \sum_{i,t} \alpha_i \Bigl( v_t^i R_t (v_t^i)^\top\Bigr),\\
\label{eq:cov-rewritten-obj-fun}J_{\text{cov}}(\covar_t^i, K_t^i)
& \Let \sum_{i,t} \alpha_i \Bigl(\operatorname{tr}(R_t K_t^i \covar_t^i (K_t^i)^\top)\Bigr)\nonumber\\
&\quad + \eps \;\GGW^2_2(\covar_N, \covar_d),
\end{align}
where the GW cost can be explicitly written as
\begin{align}
    \GGW^2_2(\covar_{\horizon},\covar_d) = 4 \bigl( \tr(\covar_{\horizon}) &- \tr(\covar_d)\bigr)^2 +  8 \norm{\covar_{\horizon}}_F ^2\nn \\& - 16 \tr(D_{\horizon}D_d). \nn
\end{align}
The diagonal matrices \(D_{\horizon}\) and \(D_d\) consists of the eigen values of \(\covar_{\horizon}\) and \(\covar_d\) arranged in a descending order.  
\subsection{Data-driven parametrization, chance constraints, and GW steering}
We now formulate our data-driven GW steering problem. Recall that our data is persistently exiting. To this end, let \(T \in \N\). For a \(T\)-duration experiment and collected state-action data \(\aset[]{x^{\dataindex}_0,\ldots,x^{\dataindex}_{T}}\) and \(\aset[]{u^{\dataindex}_0,\ldots,u^{\dataindex}_{T-1}}\) define the Hankel matrices by \(\mathcal{X}_{0;T} \Let [x^{\dataindex}_0\,\ldots\,x^{\dataindex}_{T-1}]\), \(\mathcal{U}_{0:T} \Let [u^{\dataindex}_0\,\ldots\,u^{\dataindex}_{T-1}]\) and \(\mathcal{X}_{1:T} \Let [x^{\dataindex}_1\,\ldots\,x^{\dataindex}_{T}]\).\footnote{See Appendix \ref{appen:A} for more details on data-driven parametrization.} We parametrize the feedback gain \((K_{t}^i)_{t,i}\) with \(G_t^i \in \Rbb^{T \times \stdim}\) by
\begin{align}\label{eq:gain_parametrization}
    \begin{pmatrix}
        K_t^i \\ I_n
    \end{pmatrix} = 
    \begin{pmatrix}
        \condata \\ \stdata
    \end{pmatrix} G_t^i
\end{align}
This is a standard parametrization; see \cite{ref:DDCont:main:CDP:PT,ref:DDCovSteer:JP:PT}. Our next result provides a deterministic reformulation of the chance constraints in \eqref{eq:chance_cons} under a fixed risk allocation, following \cite{ref:OnoWilliams:Risk:Allocation, ref:KumKen:CCGMM:CDC}. Based on Boole's inequality, the risk allocation decomposes the joint chance constraints into individual constraints, thereby providing tractability. We further relax the problem to handle nonconvexity by exploiting the concavity of the square root term.
\begin{myOCP}
 \begin{prop}\label{prop:chance:cons:reform}
Recall the half-plane \eqref{eq:halfplane-def} in data-\eqref{ocpdata:3}. Given a standard normal distribution, let \(\mathsf{C}(\cdot)\) be the inverse cumulative distribution function. For all \((i,j,t) \in [0;\mathsf{K}-1] \times [1;L] \times [1;\horizon-1]\) and for any \(\covar_r \in \pdmat{\stdim} \), the chance constraint in \eqref{eq:chance_cons} admits the following convex reformulation
 \begin{align}
& a_j^\top \mean_t^i + \mathsf{C}(1 - r_{ijt})\biggl(\tfrac{a_j^T \covar^i_{t}a_j}{2 \sqrt{a_j^\top \covar_r a_j}} +\tfrac{1}{2}\sqrt{a_j^\top \covar_r a_j}\biggr) \leq b_j,\label{eq:robust-cc}
    \\
    &\sum_{i=0}^{\mathsf{K}-1}
  \sum_{j=1}^{L} \sum_{t=1}^{\horizon}\alpha_i r_{ijt} \leq V,\label{eq:delta_ijk_sum}
\end{align}
where \(r_{ijt} \in \lcro{0}{1}\) is a risk parameter.
\end{prop}  
\end{myOCP}

\begin{proof}
For \(x\sim\GMM(\{\alpha_i,\mean_t^i,\covar_t^i\}_{i=0}^{\mathsf{K}-1})\), the half-plane chance constraint is implied by the set of inequalities \cite[Theorem 2]{ref:KumKen:CCGMM:CDC}
 \begin{equation}\label{eq:proof:cc}
     \begin{cases}
         a_j^T\mean_t^i + \mathsf{C}(1 - r_{ijt})\sqrt{a_j^\top\covar_t^i a_j}\le b_j,  \\
\sum_{i=0}^{\mathsf{K}-1}\sum^{L}_{j=1}\sum_{t=1}^{N} \alpha_i r_{ijt}\le V. 
     \end{cases}
 \end{equation}
To bypass the nonconvexity caused by square root term, we exploit the concavity of square root term as \cite{ref:KJ:PT:Com:Efi}, i.e., for any given \(\xi\ge 0\) and \(\xi_0\ge0\), it holds that \(\sqrt{\xi}\le\frac{1}{2\sqrt{\xi_0}}\xi + \frac{\sqrt{\xi_0}}{2}\). Therefore, the constraint is approximated by
\begin{align}\label{eq:proof:cc:concv}
    \sqrt{a_j^T \covar_{t}^i a_j} \le \frac{a_j^\top \covar^i_{t} a_j}{2 \sqrt{a_j^\top \covar_r a_j}}+\frac{\sqrt{a_j^\top \covar_r a_j}}{2}.
\end{align}
The proof follows via employing \eqref{eq:proof:cc:concv} in \eqref{eq:proof:cc}.
\end{proof}
The next theorem establishes that under the parametrization \eqref{eq:gain_parametrization} the GW-steering problem is a difference of convex program. 
\begin{myOCP}
\begin{theorem}\label{thrm:GW:DC}
Consider the OCP \eqref{eq:OCP} along with its data \eqref{ocpdata:1}--\eqref{ocpdata:3} and the parametrization \eqref{eq:gain_parametrization}. Recall the definition of the data-matrices \(\stdata\), \(\nextstdata\), and \(\condata\). Then we have the following assertion: 
\begin{enumerate}[label=\textup{(\alph*)}, leftmargin=*, widest=b, align=left]
\item \label{thrm:mean:dyn} Assume that \(\bigl(u_t^{\dataindex}\bigr)_{t=0}^{\horizon-1}\) is persistently exiciting of order \(\stdim+1\). Define the data-matrix 	
\begin{equation}
\hspace{-10mm}	\Rbb^{\stdim \times (\condim+\stdim)} \ni \conctdata \Let \nextstdata
	\begin{bmatrix}
		\condata \\
		\stdata
	\end{bmatrix}^{\dagger} = 
	\begin{bmatrix}
		\dataff & \datamu
	\end{bmatrix},
\end{equation}
where \(\datamu\in\Rbb^{\stdim\times \stdim}\) and \(\dataff\in\Rbb^{\stdim\times \condim}\). Then the mean dynamics in \eqref{eq:mean_dynamics} admits the form
\begin{align}\label{eq:meanProblem_eqConstraint2}
\datamu \mean_t^i + \dataff v_t^i - \mean^i_{t+1} = 0, 
\end{align}
for all \(t \in [0;\horizon-1]\).

\item \label{thrm:covar:dyn} For \(i \in [0;\mathsf{K}-1]\), let \(S_t^i \in \Rbb^{T \times \stdim}\) and \(Y_t^i \in \Rbb^{T \times T}\). The covariance dynamics constraint \eqref{eq:covaraicne_dynamics} can be imposed through the following equations 
\begin{align}\label{eq:covar:constraints}
\begin{cases}
\begin{pmatrix}
\covar_t^i & (S_t^i)^\top \\
S_t^i & Y_t^i
\end{pmatrix} \succeq 0,\, \covar_t^i = \stdata S_t^i,\\
\text{and }\covar_{t+1}^i = \nextstdata Y_t^i \nextstdata^\top.
\end{cases}
\end{align}
\item \label{thrm:DC:OCP}
Define the tuple \(\ol{\eta}\Let \bigl(\mean_t^i,v_t^i,\covar_t^i, Y_t^i, S_t^i\bigr)\). The GW-steering problem admits the following difference of convex programming structure
\begin{equation}\tag{(DCOCP)}
\label{eq:final:DC:OCP}
\begin{aligned} 
& \hspace{-13mm}\inf_{\ol{\eta}}	&&   \hspace{-14mm}\sum_{i,t} \alpha_i \, \bigl\{v_t^i R_t (v_t^i)^\top  + \tr \bigl(R_t \condata Y_t^i \condata^\top)\bigr\} \nn \\&& \quad & + \eps\; \GGW_2^2 (\Sigma_N,\Sigma_d)\\ 
&  \hspace{-13mm} \sbjto		&&  \hspace{-14mm}\begin{cases}
x_0 \sim \GMM(\{\alpha_i,\mu_0^i,\Sigma_0^i\}_{i=0}^{K-1}),\\
\begin{pmatrix}
        \covar_t^i & (S_t^i)^\top \\
        S_t^i & Y_t^i
    \end{pmatrix} \succeq 0,\, \covar_t^i = \stdata S_t^i , \\
    \covar_{t+1}^i = \nextstdata Y_t^i \nextstdata^\top,\,
   \text{constraints }\ref{eq:meanProblem_eqConstraint2},\\ \ref{eq:robust-cc},\,\ref{eq:delta_ijk_sum},\,
   \mean_{N}^i = \mean_{N}^j \text{ and } \covar_{N}^i = \covar_{N}^j \\ \text{ for all i,j} \in [0;\mathsf{K}-1].
\end{cases}
\end{aligned}
\end{equation}
\end{enumerate}
\end{theorem}
\end{myOCP}

\begin{proof}
The input signal \(\bigl(u_t^{\dataindex}\bigr)_{t=0}^{\horizon-1}\) is persistently exiciting of order \(\stdim+1\) and thus the system \eqref{eq:mean_dynamics} has the following equivalent representation \cite[Theorem 2]{ref:DDCovSteer:JP:PT} given by
\begin{equation}\label{eq:mean-dynamics-dd}
\mean^i_{t+1} = \nextstdata
		\begin{bmatrix}
			\condata \\
			\stdata
		\end{bmatrix}^{\dagger}
		\begin{bmatrix}
			v^i_t \\
			\mean^i_t
		\end{bmatrix},\nn
    \end{equation}
from which the first assertion readily follows.

For the second assertion, following \cite{ref:DDCovSteer:JP:PT}, we write the dynamics \eqref{eq:covariance_objective} using the parametrization \eqref{eq:gain_parametrization} for each \(i \in [0;\mathsf{K}-1]\) and \(t\in [0;\horizon-1]\), as
\begin{align}\label{eq:nonconv:cov:dyn}
    \covar_{t+1}^i = \nextstdata G_t^i \Sigma_t^i (G_t^i)^\top \nextstdata^\top.
\end{align}
We define the new decision variables
\[
S_t^i \Let G_t^i \covar_t^i \in \Rbb^{T \times \stdim}, \quad Y_t^i \succeq S_t^i (\covar_t^i)^{-1} (S_t^i)^\top,
\]
to relax the constraint \eqref{eq:nonconv:cov:dyn} into a convex one, which yields the constraints in \eqref{eq:covar:constraints} when written in the standard Schur complement form.

We show the final assertion. Notice that the first term of \eqref{eq:cov-rewritten-obj-fun} can be upper-bounded using the new decision variables:
\begin{align*}
    \tr(R_t K_t^i \covar_t^i (K_t^i)^\top) \leq \tr\bigl( R_t \condata Y_t^i \condata^\top \bigr),
\end{align*}
and consequently 
\(\ol{\eta} \mapsto \widehat{J}(\eta) \Let \sum_{i,t} \alpha_i \, \bigl\{v_t^i R_t (v_t^i)^\top  + \tr \bigl(R_t \condata Y_t^i \condata^\top)\bigr\}\)
is a convex function. Note that first two terms in 
\begin{align}
    \GGW_2^2(\covar_{\horizon},\covar_d) = 4 \bigl( \tr(\covar_{\horizon}) &- \tr(\covar_d)\bigr)^2 +  8 \norm{\covar_{\horizon}}^2_{F}\nn \\& - 16 \tr(D_{\horizon}D_d). \nn
\end{align}
are convex and the last term \(\covar_{\horizon} \mapsto \mathcal{G}(\covar_{\horizon}) \Let \tr \bigl(D_{\horizon} D_d\bigr)\) is also convex. To see this, let \(\mathcal{O}(\stdim;\Rbb) \Let \aset[]{A \in \Rbb^{\stdim \times \stdim} \suchthat AA^{\top} = I_{\stdim}}\) be the space of \(\stdim\times \stdim\) orthogonal matrices with real entries; consider the maximization problem 
\begin{equation}
    \label{eq:max_prob}
    \begin{aligned}
        & \sup_{\mathcal{S}}
        &&\tr \bigl(\mathcal{S}\covar_{\horizon}\mathcal{S}^{\top}\covar_d\bigr) \\
        & \sbjto   && \mathcal{S} \in \mathcal{O}(\stdim;\Rbb).
    \end{aligned}
\end{equation}
Note that \eqref{eq:max_prob} admits a global optimal solution; indeed, \(\mathcal{S} \mapsto \tr \bigl(\mathcal{S}\covar_{\horizon}\mathcal{S}^{\top}\covar_d\bigr)\) is continuous and \(\mathcal{O}(\stdim;\Rbb)\) is compact. Thus, existence of solution follows from \cite[Box 1.3, pp.9--10]{ref:santambrogio2023course}, and \(\covar_{\horizon} \mapsto \mathcal{G}(\covar_{\horizon}) \Let \tr \bigl(D_{\horizon} D_d\bigr)\) is the optimal value associated with the problem \eqref{eq:max_prob}. Thus \(\mathcal{G}(\cdot)\) is a maximum of linear functions (note that \(\tr(\cdot)\) is a linear operator) and thus it is convex. 

Thus the objective in \eqref{eq:final:DC:OCP} is a difference of convex function and all the constraints are convex, thus \eqref{eq:final:DC:OCP} is a difference of convex program. Our proof is complete.
\end{proof}

\begin{rem}
There are several numerical algorithms to tackle the DC program \eqref{eq:final:DC:OCP}; for example, the DC algorithm (DCA) \cite{ref:DCA:tao1997convex} (which iteratively builds a convex upper approximation of the objective function, minimizes this approximation, and then refines it using the minimizer from the previous iteration) or the Concave-Convex Procedure (CCCP) \cite{ref:CCCP:yuille2003concave} (which iteratively approximates the non-convex problem by linearizing the concave part, transforming it into a convex problem) or branch-and-bound type methods. We will specifically employ the DCA algorithm in our numerical studies because in our case the convex upper approximation is a semi-definite program (SDP). Indeed, under a proper eigen-decomposition of \(\covar_{\horizon}\) and \(\covar_d\), the term \(\mathcal{G}(\covar_N)\) admits a subgradient and consequently the convex subproblem associated with \eqref{eq:final:DC:OCP} is a SDP
\begin{equation}
\label{eq:final:DC:OCP:SDP}
\begin{aligned} 
& \inf_{\ol{\eta}}	&&   \hspace{-4mm}\widehat{J}(\eta) + 4 \bigl( \tr(\covar_{\horizon}) - \tr(\covar_d)\bigr)^2 +  8 \norm{\covar_{\horizon}} ^2_{F} ) \\&& \quad &  - 16 \tr(\covar_{\horizon}V_{N,n}^{\top}\covar_d V_{\horizon,n})\nn\\ 
&   \sbjto		&&  \hspace{-4mm}\begin{cases}
\text{constraints in \eqref{eq:final:DC:OCP} hold,}
\end{cases}
\end{aligned}
\end{equation}
where \(V_{N,n}\) is the matrix obtained by decomposing \(\covar_{\horizon}\) in the \(n\)-th iteration of the DCA algorithm; we refer \cite{ref:KM:KK:MinDenLin:2024} for more details. 
\end{rem}
\section{Numerics}\label{sec:numerics}
\label{sec:num_exp}

We illustrate our result via a numerical example. Consider the linear system 
\begin{align}\label{eq:num:dynamics}
    x_{t+1} = \begin{pmatrix}1 & 0.1 \\ -0.3 & 1.0 \end{pmatrix} x_t + \begin{pmatrix}  0.7 \\ 0.4
    \end{pmatrix} u_t,
\end{align}
for all \(t \in \Nz\). The parameters corresponding to the \(\GMM\) are fixed as follows: \(\mathsf{K}=3\) and \(\alpha = (0.4,0.3,0.3)\); the initial means are \(\mean^1_0 = (-5.0, 0.0)\), \(\mean^2_0 = (-1.0,0.0)\), and \(\mean^3_0 = (-9.0,0.0)\); for each \(i=1,2,3\), the initial covariances are \(\covar_0^i = I_2\). The covariance of the target Gaussian distribution is given by 
\[
\covar_d = \begin{pmatrix}
    2 & 0 \\ 
    0 & 0
\end{pmatrix},
\]
i.e., the terminal covariance is to be steered to a straight line. For the chance constraints in \eqref{eq:chance_cons}, we fix \(V = 0.01\) and \(r_{ijt} = \frac{V}{N}\). For the state chance constraints, we set \(a_1 = (0.707,0.707)\) and \(b_1 = 8.48\) which defines a half-plane in \(\Rbb^2\). Finally, fixing a time horizon \(\horizon = 15\) and with \(T = 10\), we consider the minimum energy optimal control problem \eqref{eq:final:DC:OCP}.

\begin{figure}[htbp]
    \centering
    \includegraphics[width=1.0\linewidth]{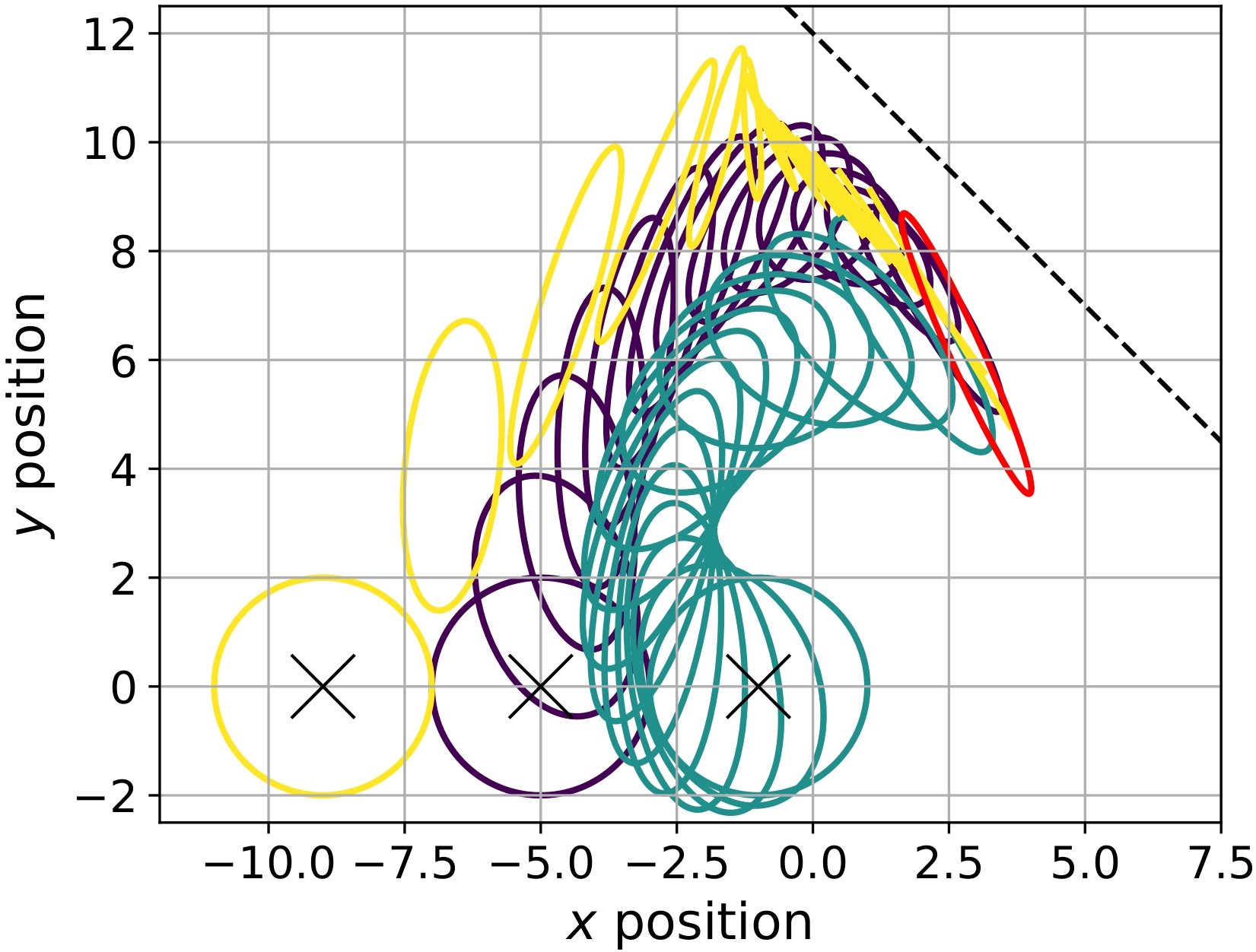}
    \caption{The trajectory of \(2\sigma\)-range of \(\covar^i_t\) obtained by our approach. The initial mean of each trajectory is denoted by a cross, and the black dotted line specifies the state constraint. The terminal \(2\sigma\) is represented by a red ellipsoid.
}
\label{fig:dd-2sigma-traj}
\end{figure}

\begin{figure}[htbp]
    \centering
    \includegraphics[width=1.0\linewidth]{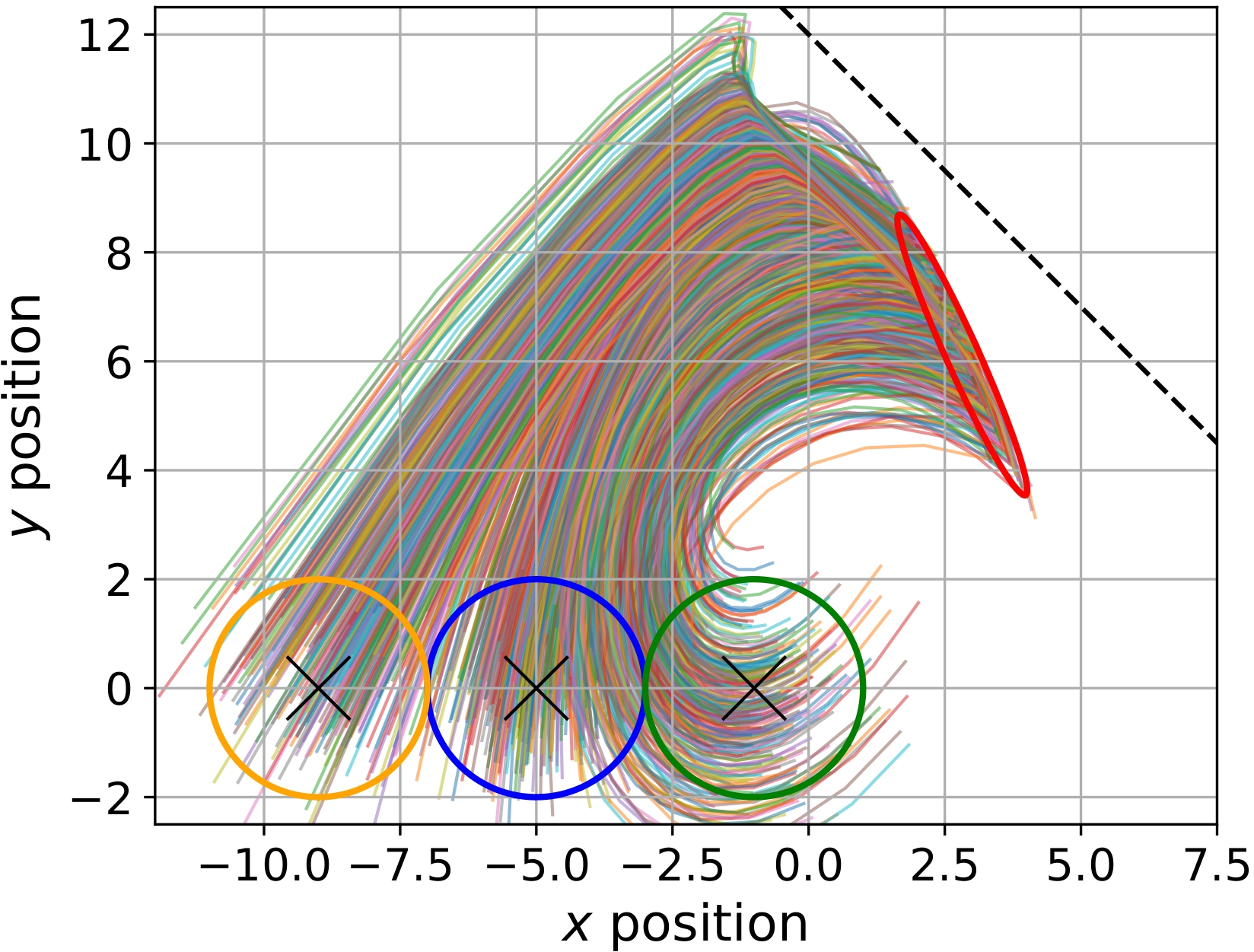}
    \caption{1000 Monte Carlo sampled trajectories obtained by our method.
}
    \label{fig:dd-sample-traj}
\end{figure}

We solved the DC program using the MOSEK solver \cite{ref:mosek} and CVXPY modeller. Figure \ref{fig:dd-2sigma-traj} depicts the trajectory of 2\(\sigma\) using the proposed data-driven approach, and Figure \ref{fig:dd-sample-traj} shows \(1000\) Monte Carlo sampled trajectories. It can be seen that the samples are successfully steered from a multi-modal initial distribution to the terminal Gaussian distribution while satisfying the chance constraint throughout the time horizon and achieves the shape of the specified straight line. We now present a numerical comparison with the data-driven covariance steering algorithm \cite{ref:DDCovSteer:JP:PT} and also examine the scenario where the data is corrupted by noise.
\subsection*{Comparison with the standard covariance steering}
\begin{figure}[htbp]
    \centering
    \includegraphics[width=1.0\linewidth]{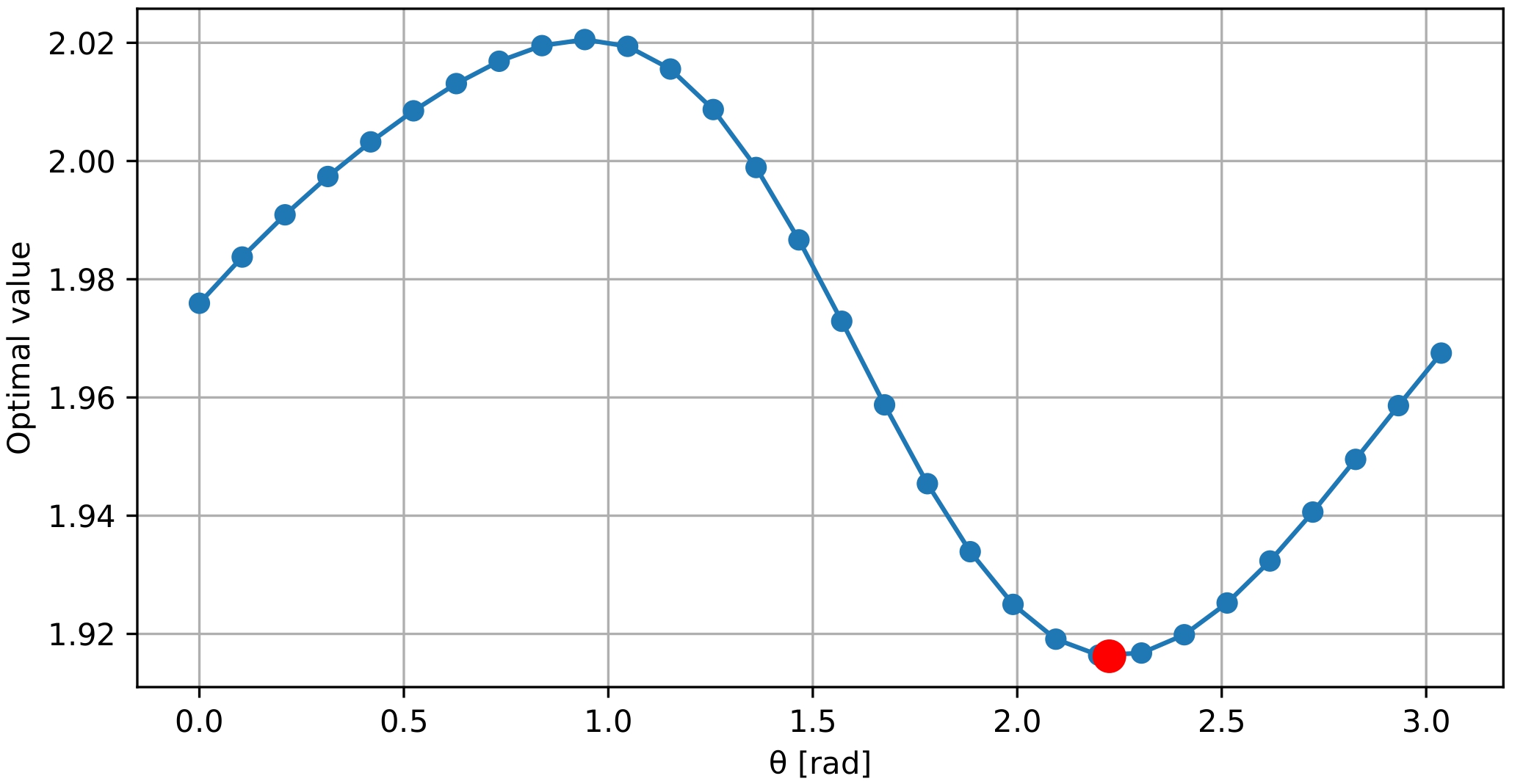}
    \caption{The blue dots stand for the optimal value of \ref{eq:standard-covariance-steering} while the red dot is the control cost achieved by the proposed approach.}
    \label{fig:opt-val-angle}
\end{figure}

To see the optimality of the proposed approach, we implement the standard covariance steering with varying angles of the terminal covariance and fixed mean. In the previous simulation without the chance constraint, the optimal terminal mean is \(\mean_{\horizon}\as = (2.320, 6.922)\). Consider the following standard covariance steering problem
\begin{equation}\label{eq:standard-covariance-steering}
\begin{aligned}
&  \hspace{-3mm}\inf_{\pi(\cdot)}	&&   \hspace{-2mm} \expect \biggl[ \sum_{t=0}^{\horizon-1} \norm{u_t}_{R_t}^2 \biggr]  \\
&  \hspace{-2mm} \sbjto		&&  \hspace{-3mm}\begin{cases}
\text{dynamics }\eqref{LTI_system},\, x_0 \sim \GMM\{\alpha_i,\mu_0^i,\Sigma_0^i\}_{i=0}^{K-1},\\
u_k = \pi_k(x_k), \text{ and } x_N \sim \mathcal{N}(\mu_N^*, \Sigma_r(\theta)),
\end{cases}
\end{aligned}
\end{equation}
where the terminal covariance is given by
\begin{align*}
    \Sigma_r(\theta) = r(\theta)^\top \Sigma_r r(\theta),\quad r(\theta) \Let \begin{pmatrix}
        \cos \theta  & -\sin \theta  \\
        \sin \theta & \cos \theta
    \end{pmatrix}.
\end{align*}
We solve this OCP using the technique given in \cite{ref:ACC:balci2024density}. Figure. \ref{fig:opt-val-angle} presents the optimal value of \ref{eq:standard-covariance-steering} for different $\theta$. It can be observed that the proposed approach successfully finds the optimal angle that achieves the minimal control cost, which presents the virtue of the invariance of GW cost with respect to rotation.
\subsection*{Regularized approach for noise corrupted data}
In practical applications, the collected data is often corrupted by noise. To attenuate this issue, some works analyze and mitigate the effect of noise in data-driven control. In contrast to the proposed direct approach, one option is to first identify the system from data and then synthesize a controller, which is referred to as the indirect design. This approach gives the exact system matrices in noiseless case as the assumption of the persistence of excitation coincides the identifiability. The simplest method for handling noisy data is to ignore the noise and apply the indirect design --- this is known as the \emph{certainty equivalence} approach. In \cite{On-the-Role-of-Regularization}, it is shown that this approach exhibits a certain degree of robustness against noise at the expense of optimality to obtain a stabilizing controller in the LQR setting. Let \(\mathcal{L} \Let \begin{pmatrix}
        \condata \\ \stdata
    \end{pmatrix}\),
and let \(\Pi_{\mathcal{L}}\Let I_T - \mathcal{L}^\dagger\mathcal{L}\) be the orthogonal projection on the nullspace of \(\mathcal{L}\) and \(I_T\) is the \(T \times T\) identity matrix. In our setting, the equivalent problem is given by the following OCP
\begin{equation}\label{eq:indirect-design}
\begin{aligned} 
&  \inf_{\ol{\eta}}	&&   \text{the cost function in \ref{eq:final:DC:OCP}}\\ 
&   \sbjto && \hspace{-1mm}\begin{cases}
\text{the constraints in \ref{eq:final:DC:OCP} }, \\
\Pi_\mathcal{L}S_k^i = 0,
\end{cases}
\end{aligned}
\end{equation}
for all \(i \in [0;\mathsf{K}-1]\) and \(t \in [0;\horizon-1]\). In \cite{ref:DDCovSteer:JP:PT}, it was recorded that direct data driven approach achieves better optimality while sacrificing the robustness, and a regularized approach, blend of certainty equivalence and direct data-driven approach, trade-off those effects, which is posed as the following problem with the hyperparameter \(\lambda \geq 0\)
\begin{equation}\label{eq:regularized-formulation}
\begin{aligned} 
& \inf_{\ol{\eta}}	&&  \text{the cost function in \ref{eq:final:DC:OCP}} + \lambda \norm{\Pi_\mathcal{L}S_t^i}\\ 
& \sbjto && \hspace{-1mm}\begin{cases}
\text{the constraints in \ref{eq:final:DC:OCP}}.
\end{cases}
\end{aligned}
\end{equation}
It can be shown that, for sufficiently large \(\lambda\), the problem coincides the indirect design \cite{On-the-Role-of-Regularization}, i.e., \ref{eq:indirect-design}. We briefly provide an intuitive explanation of the role of regularization in our context, along with the results of numerical simulations, and leave rigorous mathematical arguments for future work. Denote the noise in the data by \(\mathcal{N}\in \Rbb^{T\times\stdim}\) i.e., it holds that \(\nextstdata - \mathcal{N} = \begin{pmatrix}
    B & A
\end{pmatrix} \begin{pmatrix}
    \condata \\ \stdata
\end{pmatrix}.\)
Using the parametrization \eqref{eq:gain_parametrization}, the covariance dynamics that corresponds to \eqref{eq:nonconv:cov:dyn} is given by
\begin{align}\label{eq:noisy-cov-dynamics}
    \covar_{t+1}^i = (\nextstdata - \mathcal{N})G_t^i \Sigma_t^i (G_t^i)^\top (\nextstdata - \mathcal{N})^\top.
\end{align}
Considering the difference between the right-hand sides of \eqref{eq:nonconv:cov:dyn} and \eqref{eq:noisy-cov-dynamics}, choosing smaller values for \(G_t^i\) --- which is implied by smaller values for \(S_t^i\) in the new decision variables --- enables the noisy covariance dynamics to more closely approximate the ideal covariance dynamics \eqref{eq:nonconv:cov:dyn}. However, this comes at the cost of imposing additional constraints on the original problem, a trade-off that becomes especially prominent for large \(\lambda\). Hence, the regularized formulation \eqref{eq:regularized-formulation} is expected to provide a performance trade-off, for carefully chosen parameters. 
\begin{table}[b]
\centering
\caption{Indirect}
\label{tab:comparison-stats-indirect}
{\fontsize{10pt}{12pt}\selectfont
\begin{tblr}{l c c c}
\hline[2pt]
\SetRow{azure9}
\(\beta\) & \(\text{Control cost}\) & GGW cost & violation of CC (\%)\\
\hline
\(10^{-4}\) & 10.039 & 0.0031 & 0.0 \\
\hline
\(10^{-3}\) & 10.019 & 0.0031 & 0.0 \\
\hline
\(10^{-2}\) & 9.984 & 0.0031 & 0.0 \\
\hline[2pt]
\end{tblr}
}
\end{table}
Tables \ref{tab:comparison-stats-indirect}, \ref{tab:comparison-stats-reg}, and \ref{tab:comparison-stats-direct} present a comparison of the indirect, regularized, and direct data-driven approaches. The values in the tables are the averages of 10 independent trials, conducted under the same settings as the line alignment problem, with the hyperparameter \(\lambda = 7.0 \times 10^{-2}\) for the regularized problem. The violation of the chance constraint is computed using 1000 Monte Carlo samples for each trial. Here, \(\beta\) specifies the intensity of the noise in data; that is, a  noise \(w_t \sim \mathcal{N}(0, \beta I_{\stdim})\) is added to the dynamics during data collection.  For the all method, we obtain almost the same GGW cost at terminal step. It can be observed that the indirect approach is the most conservative, yielding the worst control cost while ensuring that the chance constraint is respected. On the other hand, the direct data-driven approach violates the chance constraint due to overfitting to the corrupted data, with the severity of the violation increasing as the noise level increases. The regularized approach, by contrast, strikes a balance between optimality and robustness against noise.

\begin{table}[t]
\centering
\caption{Regularized}
\label{tab:comparison-stats-reg}
{\fontsize{10pt}{12pt}\selectfont
\begin{tblr}{l c c c}
\hline[2pt]
\SetRow{azure9}
\(\beta\) & \(\text{Control cost}\) & GGW cost & violation of CC (\%)\\
\hline
\(10^{-4}\) & 10.021 & 0.0054 & 0.0 \\
\hline
\(10^{-3}\) & 9.828 & 0.0054 & 0.0 \\
\hline
\(10^{-2}\) & 8.794 & 0.0050 & 0.0 \\
\hline[2pt]
\end{tblr}
}
\end{table}

\begin{table}
\centering
\caption{Direct}
\label{tab:comparison-stats-direct}
{\fontsize{10pt}{12pt}\selectfont
\begin{tblr}{l c c c}
\hline[2pt]
\SetRow{azure9}
\(\beta\) & \(\text{Control cost}\) & GGW cost & violation of CC (\%)\\
\hline
\(10^{-4}\) & 9.958 & 0.0052 & 0.0 \\
\hline
\(10^{-3}\) & 9.013 & 0.0024 & 3.4 \\
\hline
\(10^{-2}\) & 8.458 & 0.0 & 16.8 \\
\hline[2pt]
\end{tblr}
}
\end{table}





\section{Conclusion}\label{sec:conclusion}
We developed a numerical algorithm for data-driven density steering using the Gromov-Wasserstein terminal cost. We demonstrated that for an unknown linear system, in the presence of chance constraints and under a specific control parametrization, the resulting optimal control problem can be formulated as a difference-of-convex program. We then solved this problem using the Difference-of-Convex Algorithm (DCA).

In future work, we aim to extend our approach to noisy and unknown linear systems, exploring density steering in both stochastic and robust settings. As mentioned in the introduction, covariance steering can handle uncertainty without heavy online computation; however, existing works highly depend on the availability of an accurate model, which is not necessarily given. The theoretically guaranteed methodology to tackle this problem remains relatively unexplored, with existing methods relying on relaxation techniques, whereas we will employ semi-infinite optimization algorithms, following \cite{ref:GanchaLCSS}, to achieve more reliable solutions.

\section{Appendix A}\label{appen:A}
\subsection{Data-Driven Parametrization}
We record some standard definitions from the data-driven control literature; we direct the readers to \cite{ref:DDCont:main:CDP:PT} for more details and application in density steering problems. Let \(\sigma \in \N\). Given a signal \(z:\Z \lra \Rbb^{\sigma}\), i.e., \(\Z \ni k \mapsto z_k\) for each \(k\), its Hankel matrix is given by
\begin{align}
        \mathcal{H}_{i,p,q}\Let 
        \begin{pmatrix} 
        z_{i} & z_{i+1} & \dots & z_{i+q-1} \\
        z_{i+1} & z_{i+2} & \dots & z_{i+q} \\
        \vdots & \vdots & \ddots & \vdots \\
        z_{i+p-1} & z_{i+p} & \dots &  z_{i+p+q-2} \\
    \end{pmatrix}.
    \end{align}
where \((i,p,q) \in \Z \times \Nz \times \Nz\). Note that for \(p=1\) the matrix \(\mathcal{H}_{i,p,q}\) consists of only one block row, in that case we write \(\mathcal{H}_{i,1,q} \Let \mathcal{H}_{i,p} = [z_i\,z_{i+1}\,\ldots\,z_{i+p-1}]\).

\begin{defn}
 \textbf{(Persistence of Excitation)} The signal \(\lcrc{0}{T-1}\cap \Z \ni k \mapsto z_k \in \Rbb^{\sigma}\) is \emph{persistently exciting} of order \(\ell\) if \(\mathsf{rank}\bigl(\mathcal{H}_{0,\ell,T-\ell+1} \bigr) = \sigma \ell\).
\end{defn}
It is well known \cite{ref:DDCont:main:CDP:PT} that in order for a signal to be persistently exciting of order \(\ell\), it must be sufficiently long, i.e., it must hold that \(T \geq (\sigma + 1)\ell -1\). The following result is the centerpiece in data-driven control literature.
\begin{lemm} \label{lemm:willem_funda}
\textbf{(Willem's Fundamental Lemma)}
Consider the discrete-time system \(x_{t+1}= Ax_t+B u_t\) and let the pair \((A,B)\) be controllable. For a \(T\)-duration experiment and collected state-action data \(\aset[]{x^{\dataindex}_0,\ldots,x^{\dataindex}_{T}}\) and \(\aset[]{u^{\dataindex}_0,\ldots,u^{\dataindex}_{T-1}}\) define the Hankel matrices (with \(p=1\)) by \(\mathcal{X}_{0;T} \Let [x^{\dataindex}_0\,\ldots\,x^{\dataindex}_{T-1}]\), \(\mathcal{U}_{0:T} \Let [u^{\dataindex}_0\,\ldots\,u^{\dataindex}_{T-1}]\) and \(\mathcal{X}_{1:T} \Let [x^{\dataindex}_1\,\ldots\,x^{\dataindex}_{T}]\). If the control signal \(\{u_t\}_{t=0}^{T-1}\) is persistently exiting of order \(n+1\), then
\begin{align}
 \textsf{rank} \begin{pmatrix}
            \mathcal{U}_{0:T} \\ \mathcal{X}_{0:T}
        \end{pmatrix} = n+m.\nn
    \end{align}
\end{lemm}
\bibliographystyle{ieeetr}
\bibliography{references}

\begin{thebibliography}{10}

\bibitem{ref:DDCont:main:CDP:PT}
C.~D. Persis and P.~Tesi, ``Formulas for data-driven control: Stabilization, optimality, and robustness,'' {\em IEEE Transactions on Automatic Control}, vol.~65, no.~3, pp.~909--924, 2020.

\bibitem{ref:DZ:JR:PT:AA:BSP:21}
D.~Zheng, J.~Ridderhof, P.~Tsiotras, and A.~Agha-mohammadi, ``Belief space planning: A covariance steering approach,'' in {\em 2022 International Conference on Robotics and Automation (ICRA)}, pp.~11051--11057, IEEE, 2022.

\bibitem{ref:JK:PT:CovUnknParam:2023}
J.~W. Knaup and P.~Tsiotras, ``Covariance steering for systems subject to unknown parameters,'' in {\em 2023 62nd IEEE Conference on Decision and Control (CDC)}, pp.~1790--1795, 2023.

\bibitem{ref:GR:PT:StocCont:UAV}
G.~Rapakoulias and P.~Tsiotras, ``Stochastic control of uavs: An optimal tradeoff between performance, flight smoothness and control effort,'' 2024.
\newblock doi: \url{https://arxiv.org/abs/2409.10369}.

\bibitem{ref:DDCovSteer:JP:PT}
J.~Pilipovsky and P.~Tsiotras, ``Data-driven covariance steering control design,'' in {\em 2023 62nd IEEE Conference on Decision and Control (CDC)}, pp.~2610--2615, IEEE, 2023.

\bibitem{ref:L4DC:DDRobCovSteer:JP:PT}
J.~Pilipovsky and P.~Tsiotras, ``Data-driven robust covariance control for uncertain linear systems,'' in {\em Proceedings of the 6th Annual Learning for Dynamics \& Control Conference}, vol.~242 of {\em Proceedings of Machine Learning Research}, pp.~667--678, PMLR, 15--17 Jul 2024.

\bibitem{ref:KM:KK:MinDenLin:2024}
K.~Morimoto and K.~Kashima, ``Minimum energy density steering of linear systems with {G}romov--{W}asserstein terminal cost,'' {\em IEEE Control Systems Letters}, vol.~8, pp.~586--591, 2024.

\bibitem{ref:NH:SG:KM:KK:Formation:shape:GW}
H.~Nakashima, S.~Ganguly, K.~Morimoto, and K.~Kashima, ``Formation shape control using the gromov-wasserstein metric,'' 2025.
\newblock doi: \url{https://arxiv.org/abs/2503.21538}.

\bibitem{ref:IMB:BE:GMM:2024}
I.~M. Balci and E.~Bakolas, ``Density steering of gaussian mixture models for discrete-time linear systems,'' in {\em 2024 American Control Conference (ACC)}, pp.~3935--3940, 2024.

\bibitem{ref:KumKen:CCGMM:CDC}
N.~Kumagai and K.~Oguri, ``Chance-constrained gaussian mixture steering to a terminal gaussian distribution,'' in {\em 2024 IEEE 63rd Conference on Decision and Control (CDC)}, pp.~2207--2212, 2024.

\bibitem{ref:OnoWilliams:Risk:Allocation}
M.~Ono and B.~C. Williams, ``Iterative risk allocation: A new approach to robust model predictive control with a joint chance constraint,'' in {\em 2008 47th IEEE Conference on Decision and Control}, pp.~3427--3432, 2008.

\bibitem{ref:GGW:delon2022gromov}
J.~Delon, A.~Desolneux, and A.~Salmona, ``Gromov--{W}asserstein distances between {G}aussian distributions,'' {\em Journal of Applied Probability}, vol.~59, no.~4, pp.~1178--1198, 2022.

\bibitem{ref:FM:GW:2011}
F.~M{\'e}moli, ``Gromov-{W}asserstein distances and the metric approach to object matching,'' {\em Foundations of {C}omputational {M}athematics}, vol.~11, pp.~417--487, 2011.

\bibitem{ref:ACC:balci2024density}
I.~M. Balci and E.~E.~Bakolas, ``Density steering of {G}aussian mixture models for discrete-time linear systems,'' in {\em 2024 American Control Conference (ACC)}, pp.~3935--3940, IEEE, 2024.

\bibitem{ref:KJ:PT:Com:Efi}
J.~W. Knaup and P.~Tsiotras, ``Computationally efficient covariance steering for systems subject to parametric disturbances and chance constraints,'' in {\em 2023 62nd IEEE Conference on Decision and Control (CDC)}, pp.~1796--1801, 2023.

\bibitem{ref:santambrogio2023course}
F.~Santambrogio, {\em A {C}ourse in the {C}alculus of {V}ariations: {O}ptimization, {R}egularity, and {M}odeling}.
\newblock Universitext, Springer Nature, 2023.

\bibitem{ref:DCA:tao1997convex}
P.~D. Tao and L.~H. An, ``Convex analysis approach to {DC} programming: theory, algorithms and applications,'' {\em Acta mathematica vietnamica}, vol.~22, no.~1, pp.~289--355, 1997.

\bibitem{ref:CCCP:yuille2003concave}
A.~L. Yuille and A.~Rangarajan, ``The concave-convex procedure,'' {\em Neural computation}, vol.~15, no.~4, pp.~915--936, 2003.

\bibitem{ref:mosek}
{\relax MOSEK ApS}, ``The {MOSEK} optimization toolbox for {MATLAB} manual. version 9.0.,'' 2019.
\newblock url: \url{http://docs.mosek.com/9.0/toolbox/index.html}.

\bibitem{On-the-Role-of-Regularization}
F.~D\"{o}rfler, P.~Tesi, and C.~D. Persis, ``On the role of regularization in direct data-driven {LQR} control,'' in {\em 2022 IEEE 61st Conference on Decision and Control (CDC)}, pp.~1091--1098, 2022.

\bibitem{ref:GanchaLCSS}
S.~Ganguly and D.~Chatterjee, ``Exact solutions to minmax optimal control problems for constrained noisy linear systems,'' {\em IEEE Control Systems Letters}, vol.~8, pp.~2063--2068, 2024.

\end{thebibliography}

\end{document}